# The Long Term Fréchet distribution: Estimation, Properties and its Application

Pedro Luiz Ramos, Diego Nascimento and Francisco Louzada

Institute of Mathematical Science and Computing, University of São Paulo,
São Carlos, Brazil

**Abstract**

In this paper a new long-term survival distribution is proposed. The so called long term Fréchet distribution allows us to fit data where a part of the population is not susceptible to the event of interest. This model may be used, for example, in clinical studies where a portion of the population can be cured during a treatment. It is shown an account of mathematical properties of the new distribution such as its moments and survival properties. As well is presented the maximum likelihood estimators (MLEs) for the parameters. A numerical simulation is carried out in order to verify the performance of the MLEs. Finally, an important application related to the leukemia free-survival times for transplant patients are discussed to illustrates our proposed distribution

## 1. Introduction

Extreme value models play an important role in statistic. The generalized extreme value (GEV) distribution (Jenkinson, 1955) and its sub-models are widely used in application involving extreme events. These sub-models are the well known Weibull, Fréchet and Gumbel distributions. The Fréchet distribution can be seen as the inverse Weibull distribution which gives a probability density function (PDF) such as

$$f(t,\lambda,\alpha) = \frac{\alpha}{\lambda}\left(\frac{t}{\lambda}\right)^{-(\alpha+1)} e^{-\left(\frac{t}{\lambda}\right)^{-\alpha}}. \qquad (1)$$

The survival function is given by

$$S(t,\lambda,\alpha) = 1 - e^{-\left(\frac{t}{\lambda}\right)^{-\alpha}}. \qquad (2)$$

Although the GEV distribution is the most used generalization of the Fréchet model, other distributions has been proposed in the literature. De Gusmão (2012) proposed a three parameter generalized inverse Weibull distribution in which includes the Fréchet distribution. Krishna et al. (2013) proposed the Marshall-Olkin Fréchet distribution. Barreto-Souza et al. (2011) discussed some results for beta Fréchet distribution. However, in survival studies habitually the researches may consider a portion of the population as cured during a given treatment, this type of distribution is called long-term (LT) survival models.

In this study, a long-term survival novel proposing a mixture model introduced by Berkson and Gage (1952), hereafter we shall call it the long-term Fréchet distribution or simplistically the LF distribution. Some mathematical properties about the LF distribution were provided such as moments, survival properties and hazard function. The maximum likelihood estimators of the parameters and its asymptotic properties are discussed likewise. Similar studies were presented by Roman et al. (2012) for the geometric exponential distribution and by Louzada and Ramos (2017) for the weighted Lindley distribution. It was performed a numerical simulation towards to examine the performance of the MLEs. Finally, our proposed methodology is illustrated in a real data set related to the leukemia free-survival times (in years) for the 50 autologous transplant patients.

The paper is organized as follows. Section 2 presents the long term Fréchet distribution and its mathematical properties. Section 3 discusses the parameter estimation under the maximum likelihood approach. Section 4 presents a simulation study under different values of the parameters and different levels of censorship. The proposed methodology is also fully illustrated in a real data set. Lastly, Section 6 summarizes the founds in this study and its potential contribution.

## 2. Long Term Fréchet distribution

Long-term survivors are an important feature to incorporate in the modeling process, since a portion of the population may no longer be eligible to the event of interest (according to Maller and Zhou, 1995; or Perdona and Louzada, 2011). Hence the population can be segregate as a not eligible to the event of interest with probability $p$ and as eligible (in risk) to the event of interest with probability $(1 − p)$. The long-term survivor is expressed as

$$S(t; p, \boldsymbol{\theta}) = p + (1 − p)S_0(t; \boldsymbol{\theta}), \quad (3)$$

where $p \in (0,1)$ and $S_0(t; \boldsymbol{\theta})$ is the survival function related to the eligible group. The obtained survival function (not conditional) is improper and its limit corresponds to the individual proportion cure. From the survival function one can easily derive the PDF (improper) given by

$$f(t; p, \boldsymbol{\theta}) = -\frac{\partial}{\partial t} S(t; p, \boldsymbol{\theta}) = (1 − \pi)f_0(t; \boldsymbol{\theta}), \quad (4)$$

where $f_0(t; \boldsymbol{\theta})$ is the PDF related to the susceptible group.

Considering that $f_0(t; \boldsymbol{\theta})$ follows a Fréchet distribution, then the PDF of the Long Term Fréchet (LF) distribution is given by

$$f(t; \lambda, \alpha, p) = \frac{\alpha(1-p)}{\lambda}\left(\frac{t}{\lambda}\right)^{-(\alpha+1)} e^{-\left(\frac{t}{\lambda}\right)^{-\alpha}}, \quad (5)$$

where $\lambda > 0$, $\alpha > 0$ and $p \in (0,1)$. The cumulative distribution function is given by

$$F(t; \lambda, \alpha, p) = (1 − p)e^{-\left(\frac{t}{\lambda}\right)^{-\alpha}}. \quad (6)$$

In this case, the LF has the quantile function in closed-form and is given by

$$t_u = \lambda \log\left(\frac{(1 − p)}{u}\right)^{-\frac{1}{\alpha}}, \quad (7)$$

where $0 \leq u < 1$. The $r$-th moments of T about the origin is

$$E(T^r; \lambda, \alpha, p) = (1 − p)\lambda^r \Gamma\left(1 − \frac{r}{\alpha}\right), \quad \alpha > r, \quad (8)$$

for $r \in \mathbb{N}$ and $\Gamma(x) = \int_0^\infty e^{-y} y^{x-1} dy$ is called gamma function. Along with some algebraic manipulation the mean and variance of the LF distribution are given, respectively, by

$$E(T; \lambda, \alpha, p) = (1 − p)\lambda \Gamma\left(1 − \frac{1}{\alpha}\right), \quad \alpha > 1$$

and

$$V(T;\lambda,\alpha,p) = (1-p)\lambda^2\left(\Gamma\left(1-\frac{2}{\alpha}\right) - (1-p)\Gamma\left(1-\frac{1}{\alpha}\right)^2\right), \quad \alpha > 2.$$

The survival and hazard functions of $LF(\lambda,\alpha,p)$ distribution is given by

$$S(t;\lambda,\alpha,p) = p + (1-p)\left(1 - e^{-\left(\frac{t}{\lambda}\right)^{-\alpha}}\right) \quad (9)$$

and

$$h(t;\lambda,\alpha,p) = \frac{\frac{\alpha}{\lambda}\left(\frac{t}{\lambda}\right)^{-(\alpha+1)} e^{-\left(\frac{t}{\lambda}\right)^{-\alpha}}}{p + (1-p)\left(1 - \exp\left(-\left(\frac{t}{\lambda}\right)^{-\alpha}\right)\right)}. \quad (10)$$

Figure 1 shows some cases about the PDF and the survival function shapes applied to LF distribution.

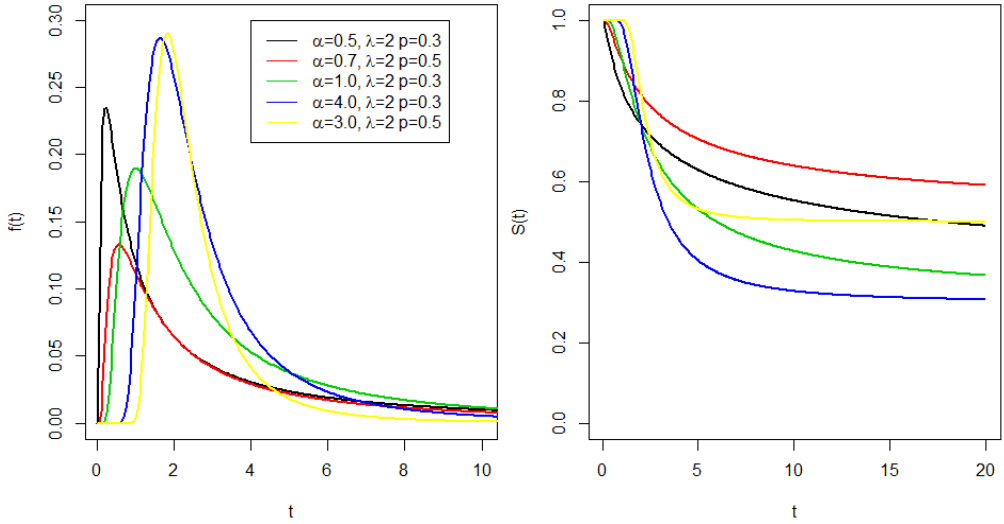

Figure 2. Left panel: Probability density function of the LF distribution. Right panel: Survival function of the LF distribution

## 3. Parameter Estimation

For each failure time related to the i-th individual, it may not be perceived or subject by the right censoring. Furthermore, the random censoring times $C_i$s are independent of $T_i$s (non-censored time) and their distribution does not depend on the parameters. In a scenario of a $n$ sample of size, the data set will be describe by $\mathcal{D} = (t_i, \delta_i)$, where $t_i = \min(T_i, C_i)$ and $\delta_i = I(T_i \leq C_i)$. This general random censoring scheme has as special case type I and II censoring mechanism. The likelihood function is given by

$$L(\boldsymbol{\theta}; \mathcal{D}) = \prod_{i=1}^{n} f(t_i; \boldsymbol{\theta})^{\delta_i} S(t_i; \boldsymbol{\theta})^{1-\delta_i}.$$

Let $T_1, \ldots, T_n$ be a random sample of LF distribution, the likelihood function considering data with random censoring is given by

$$L(\lambda, \alpha, p; \mathcal{D}) = \frac{\alpha^d (1-p)^d}{\lambda^{-d\alpha}} \prod_{i=1}^{n} t_i^{-\delta_i(\alpha+1)} \exp\left(-\sum_{i=1}^{n} \delta_i \left(\frac{t_i}{\lambda}\right)^{-\alpha}\right)$$

$$\times \left(p + (1-p)\left(1 - e^{-\left(\frac{t}{\lambda}\right)^{-\alpha}}\right)\right)^{1-\delta_i},$$

where $d = \sum_{i=1}^{n} \delta_i$. The log-likelihood function is given as

$$l(\lambda, \alpha, p; \mathcal{D}) = d \log(\alpha) + d \log(1-p) + d\alpha \log \lambda - (\alpha + 1) \sum_{i=1}^{n} \delta_i \log t_i$$

$$- \sum_{i=1}^{n} \delta_i \left(\frac{\lambda}{t_i}\right)^{\alpha} + \sum_{i=1}^{n} (1 - \delta_i) \log\left(p + (1-p)\left(1 - e^{-\left(\frac{t}{\lambda}\right)^{-\alpha}}\right)\right). \quad (11)$$

The maximum likelihood estimators (MLEs) are widely explored as statistical inferential methodology due its many desirable properties, in which includes consistency, asymptotic efficiency and invariance. The MLEs are obtained from the maximization of the log-likelihood function (11). Before we derive the MLEs of the LF, let us define the following function

$$\eta_j(\lambda, \alpha, p; \mathcal{D}) = \sum_{i=1}^{n} (1 - \delta_i) \frac{\log S(t_i; \boldsymbol{\theta})}{\partial \theta_j}, \quad j = 1, 2, 3.$$

Then, the likelihood equations are given by

$$\frac{d\alpha}{\lambda} - \alpha \lambda^{\alpha-1} \sum_{i=1}^{n} \delta_i t_i^{-\alpha} + \eta_1(\lambda, \alpha, p; \mathcal{D}) = 0,$$

$$\frac{d}{\alpha} + d \log \lambda - \sum_{i=1}^{n} \delta_i \left(\frac{\lambda}{t_i}\right)^{\alpha} \log\left(\frac{\lambda}{t_i}\right) - \sum_{i=1}^{n} \delta_i \log t_i + \eta_2(\lambda, \alpha, p; \mathcal{D}) = 0$$

and

$$\frac{d}{p-1} + \eta_3(\lambda, \alpha, p; \mathcal{D}) = 0.$$

The maximization of the log-likelihood function can be performed directly by using existing statistical packages. Further information about the numerical procedures will be discussed in the next section.

According to Migon et al. (2014), under mild conditions the obtained estimators are consistent and efficient with an asymptotically normal joint distribution given by

$$(\hat{\lambda}, \hat{\alpha}, \hat{p}) \sim N_3\big((\lambda, \alpha, p), H^{-1}(\lambda, \alpha, p)\big),$$

where $H(\lambda, \alpha, p)$ is the 3×3 observed Fisher information matrix and $H_{ij}(\lambda, \alpha, p)$ is the Fisher information given by

$$H_{ij}(\boldsymbol{\theta}) = -\frac{\partial}{\partial \theta_i \partial \theta_j} l(\boldsymbol{\theta}; \mathcal{D}), \quad i,j = 1,2,3.$$

Note that, the observed Fisher information matrix was used since it is not possible to compute the expected Fisher information matrix due its lack of closed form expression. For large samples, confidence intervals approximation can be constructed for the individual parameters $\boldsymbol{\theta_i}$ i=1,2,3, assuming a confidence coefficient $100(1-\gamma)\%$ the marginal distributions are given by

$$\hat{\theta}_i \sim N\left(\theta_i, H_{ii}^{-1}(\boldsymbol{\theta})\right), \quad i = 1,2,3.$$

## 4. Simulation Study

The maximum likelihood method efficiency was analyzed through a simulation study on the LF distribution. This procedure was conducted by computing the mean relative errors (MRE) and the mean square errors (MSE) given by

$$MRE_i = \frac{1}{N}\sum_{j=1}^{N} \frac{\hat{\theta}_{i,j}}{\theta_i}, \quad MSE_i = \frac{1}{N}\sum_{j=1}^{N} (\hat{\theta}_{i,j} - \theta_i)^2, \quad \text{for } i = 1,2,3.$$

as N is the number of estimates obtained through the MLE approach. The 95% coverage probabilities of the asymptotic confidence intervals were also evaluated. The adopted approach prioritize that the expected MLEs returns the MREs closer to one with smaller MSEs. Additionally, by considering a 95% confidence level, the interval covers the true values of θ closer to 95%. Considering scenarios with sample sizes n=(10,25,50,100, 200) and N=100,000 for the simulation study, two situations are presented by considering the proportion of cure in the population of 0.3 and 0.5. In these cases, the censored proportions are observed in different levels.

In pursuance to find the maximization of the log-likelihood function, described in the equation (11), the package called maxLik available in R developed by Henningsen and Toomet (2011) was used. The numerical results are well-behaved since was not found numerical problems using the SANN method (Simulated-annealing), such as failure evidence of convergence or end on multiple maxima. The programs can be obtained, upon request.

Table 1: MREs, MSEs, $C_{95\%}$ estimates for 100,000 considering n = (25, 50, 100, 200, 300) and 45.7% of censorship.

| $\theta$ | $\alpha = 0.5$ | | | $\lambda = 2.0$ | | | $p = 0.3$ | | | 0.457 |
|---|---|---|---|---|---|---|---|---|---|---|
| $n$ | MRE | MSE | $C_{95\%}$ | MRE | MSE | $C_{95\%}$ | MRE | MSE | $C_{95\%}$ | $M_p$ |
| 25 | 1.265 | 0.060 | 0.948 | 1.100 | 3.344 | 0.810 | 1.100 | 0.021 | 0.925 | 0.461 |
| 50 | 1.114 | 0.020 | 0.952 | 1.117 | 2.002 | 0.860 | 1.024 | 0.013 | 0.938 | 0.458 |
| 100 | 1.048 | 0.008 | 0.952 | 1.098 | 1.087 | 0.893 | 0.992 | 0.008 | 0.948 | 0.457 |
| 200 | 1.022 | 0.004 | 0.953 | 1.059 | 0.488 | 0.919 | 0.991 | 0.004 | 0.952 | 0.457 |
| 300 | 1.014 | 0.003 | 0.951 | 1.039 | 0.293 | 0.928 | 0.993 | 0.003 | 0.952 | 0.457 |

Table 2: MREs, MSEs, $C_{95\%}$ estimates for 100,000 considering n = (25, 50, 100, 200, 300) and 61% of censorship.

| | $\alpha = 0.5$ | | | $\lambda = 2.0$ | | | $p = 0.5$ | | | 0.612 |
|---|---|---|---|---|---|---|---|---|---|---|
| n | MRE | MSE | $C_{95\%}$ | MRE | MSE | $C_{95\%}$ | MRE | MSE | $C_{95\%}$ | $M_p$ |
| 25 | 1.384 | 0.144 | 0.939 | 1.263 | 8.714 | 0.788 | 1.027 | 0.022 | 0.914 | 0.612 |
| 50 | 1.158 | 0.033 | 0.940 | 1.238 | 5.967 | 0.838 | 1.001 | 0.014 | 0.928 | 0.612 |
| 100 | 1.070 | 0.013 | 0.946 | 1.156 | 2.683 | 0.877 | 0.994 | 0.007 | 0.940 | 0.612 |
| 200 | 1.031 | 0.006 | 0.951 | 1.085 | 0.841 | 0.906 | 0.994 | 0.004 | 0.948 | 0.612 |
| 300 | 1.020 | 0.004 | 0.951 | 1.056 | 0.459 | 0.921 | 0.996 | 0.002 | 0.951 | 0.612 |

Table 3: MREs, MSEs, $C_{95\%}$ estimates for 100,000 considering n = (25, 50, 100, 200, 300) and 35% of censorship.

|     | $\alpha = 2.0$ | | | $\lambda = 4.0$ | | | $p = 0.3$ | | | 0.35 |
| --- | --- | --- | --- | --- | --- | --- | --- | --- | --- | --- |
| n | MRE | MSE | $C_{95\%}$ | MRE | MSE | $C_{95\%}$ | MRE | MSE | $C_{95\%}$ | $M_p$ |
| 25 | 1.103 | 0.316 | 0.951 | 1.022 | 0.347 | 0.921 | 0.997 | 0.010 | 0.927 | 0.349 |
| 50 | 1.047 | 0.116 | 0.951 | 1.011 | 0.155 | 0.938 | 0.998 | 0.005 | 0.937 | 0.348 |
| 100 | 1.023 | 0.050 | 0.950 | 1.005 | 0.074 | 0.945 | 1.000 | 0.002 | 0.945 | 0.349 |
| 200 | 1.011 | 0.023 | 0.950 | 1.003 | 0.036 | 0.947 | 1.000 | 0.001 | 0.947 | 0.349 |
| 300 | 1.008 | 0.015 | 0.951 | 1.002 | 0.024 | 0.947 | 0.999 | 0.001 | 0.947 | 0.348 |

Table 4: MREs, MSEs, $C_{95\%}$ estimates for 100,000 considering n = (25, 50, 100, 200, 300) and 53.5% of censorship.

|     | $\alpha = 2.0$ | | | $\lambda = 4.0$ | | | $p = 0.3$ | | | 0.535 |
| --- | --- | --- | --- | --- | --- | --- | --- | --- | --- | --- |
| n | MRE | MSE | $C_{95\%}$ | MRE | MSE | $C_{95\%}$ | MRE | MSE | $C_{95\%}$ | $M_p$ |
| 25 | 1.158 | 0.619 | 0.953 | 1.034 | 0.546 | 0.910 | 0.998 | 0.011 | 0.933 | 0.535 |
| 50 | 1.068 | 0.182 | 0.953 | 1.016 | 0.230 | 0.933 | 0.999 | 0.006 | 0.942 | 0.535 |
| 100 | 1.033 | 0.075 | 0.950 | 1.007 | 0.105 | 0.942 | 1.000 | 0.003 | 0.947 | 0.535 |
| 200 | 1.016 | 0.034 | 0.950 | 1.004 | 0.051 | 0.946 | 0.999 | 0.001 | 0.947 | 0.534 |
| 300 | 1.011 | 0.022 | 0.949 | 1.002 | 0.034 | 0.948 | 1.000 | 0.001 | 0.948 | 0.535 |

The estimates obtained from Tables 1-4 for $\alpha, \lambda$ and $p$ are asymptotically unbiased, implying that MREs tend to one when n increases and the MSEs decrease to zero for n large. Analyzing the MLEs performance, with a coverage probabilities tending to 0.95, good coverage properties may be deliberated for the parameter estimators. In practical applications, those estimation procedures will be relevant as shown in the next section.

## 5 Application

In this section, we considered the data set presented by Kersey et al. (1987). The results were collected in a group of 46 patients, per years, upon the recurrence of leukemia whom received autologous marrow. Table 5 shows the full data set (+ indicates censored observations).

Table 5: Leukemia free-survival times (in years) for the 46 autologous transplant patients (where + indicates censored observations).

| 0.0301 | 0.0384 | 0.0630 | 0.0849 | 0.0877 | 0.0959 | 0.1397 | 0.1616 |
| --- | --- | --- | --- | --- | --- | --- | --- |
| 0.1699 | 0.2137 | 0.2137 | 0.2164 | 0.2384 | 0.2712 | 0.2740 | 0.3863 |
| 0.4384 | 0.4548 | 0.5918 | 0.6000 | 0.6438 | 0.6849 | 0.7397 | 0.8575 |
| 0.9096 | 0.9644 | 1.0082 | 1.2822 | 1.3452 | 1.4000 | 1.5260 | 1.7205+ |
| 1.9890+ | 2.2438 | 2.5068+ | 2.6466+ | 3.0384 | 3.1726+ | 3.4411 | 4.4219+ |
| 4.4356+ | 4.5863+ | 4.6904+ | 4.7808+ | 4.9863+ | 5.0000+ | | |

The proposed model is compared with some usual long-term survival models, such as the LT Weibull and LT weighted Lindley (Louzada and Ramos, 2017). Different discrimination criterion methods are considered: the negative of the maximum value of the likelihood function $l(\hat{\boldsymbol{\theta}}; \boldsymbol{t})$, the Akaike information criterion ($AIC = -2l(\hat{\boldsymbol{\theta}}; \boldsymbol{t}) + 2k$) and the corrected AIC ($AIC + 2k(k+1)/(n-k-1)$), where k is the number of parameters to be fitted. The best model is the one which provides the minimum criterion method values. Figure 2 presents the empirical survival function adjusted by the Kaplan-Meier estimator and different LT survival distributions.

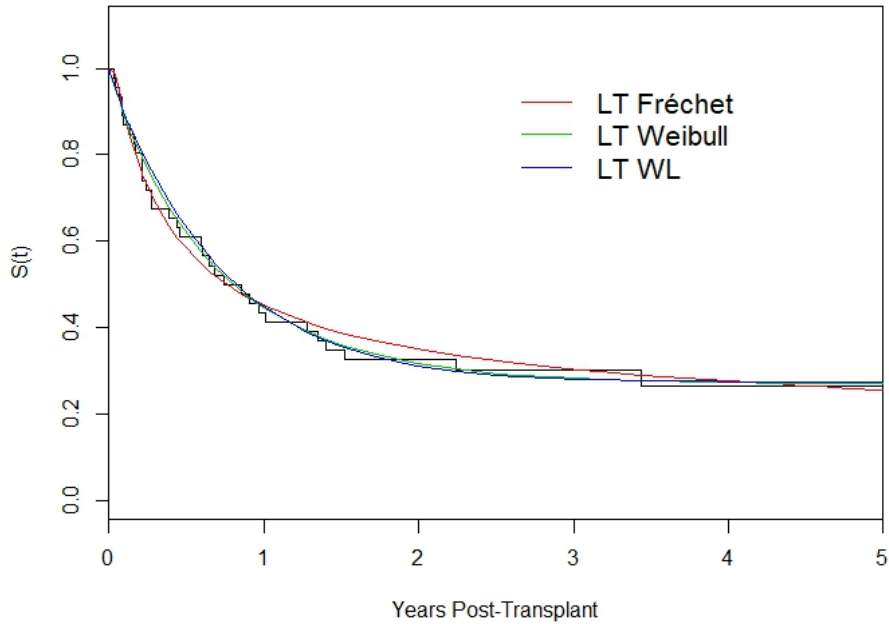

Figure 2: Survival function adjusted by the empirical survival function (Kaplan-Meier estimator), LT Fréchet, LT Weibull and LT WL distribution.

Table 6 presents the results of the different discrimination criteria for different probability distributions. Comparing the results of the different discrimination methods, we observed that the LT Fréchet distribution has better fit then the LT models under the Weibull and weighted Lindley baseline distribution.

| Method | LT Fréchet | LT Weibull | LT WL |
|---|---|---|---|
| $-\log L$ | **45.33** | 46.15 | 46.56 |
| AIC | **96.66** | 98.30 | 99.12 |
| AICc | **97.23** | 98.87 | 99.69 |

The MLEs were obtained through the same procedure as described in Section 3. The standard error (SE) and the confidence intervals, considering a 95% confidence level for $\alpha, \lambda$ and $p$ are displays in Table 7.

Table 7: MLE, Standard Error (SE), and confidence interval under 95% confidence level for $\alpha, \lambda$ and $p$.

| $\theta$ | MLE | SE | $CI_{95\%}(\theta)$ |
|---|---|---|---|
| $\alpha$ | 0.65682 | 0.01975 | (0.38140; 0.93225) |
| $\lambda$ | 0.31358 | 0.01531 | (0.07106; 0.55609) |
| $p$ | 0.12476 | 0.01597 | (0.00000; 0.37245) |

Note that, in Kersey et al. (1987) they use the non-parametric KM estimate of the cure fraction in which was $0.20$ where $(0.08; 0.32)$ is the 95% confidence interval. Therefore, results showed to be consistence with Kersey et al. results while our estimate was contained in the non-parametric interval. By using our parametric model the estimate obtained for $p$ was $0.125$ showing an overestimation of the long term survival patients. As it can be seen, through our proposed methodology the data related to the leukemia free-survival times (in months) for the 50 autologous transplant patients can be described by the LF distribution.

## 6 Discussion

In this paper, we have proposed a new long-term survival distribution called long term Fréchet distribution and its mathematical properties were studied. It was presented results towards the maximum likelihood parameters' estimators and their asymptotic properties. The estimators' efficient were present in the simulation study as the MLEs for the three unknown parameters obtained acceptable results even for small sample sizes. As such of the real dataset problem, related to the leukemia, free-survival times (in months) for the 50 autologous transplant patients. Many extensions from this present work can be considered, for instance, the parameters estimation may also be studied under an objective Bayesian analysis (Ramos et al., 2014, 2017) or using different classical methods (Louzada et al., 2016; Bakouch et al. 2017). Other approach should be to include covariates under the assumption of Cox model, i.e., proportional hazards. In conclusion, this regression model can be extended for the Bayesian approach as well.

## Acknowledgements

The authors are thankful to the Editorial Board and to the reviewers for their valuable comments and suggestions which led to this improved version.

## References

Bakouch, H. S., Dey, S., Ramos, P. L., & Louzada, F. (2017). Binomial-exponential 2 Distribution: Different Estimation Methods and Weather Applications. Trends in Applied and Computational Mathematics, 18(2), 233.

Barreto-Souza, W., Cordeiro, G. M., &Simas, A. B. (2011). Some results for beta Fréchet distribution. Communications in Statistics-Theory and Methods, 40(5), 798-811.

Berkson, J., & Gage, R. P. (1952). Survival curve for cancer patients following treatment.Journal of the American Statistical Association, 47(259), 501-515.

De Gusmão, Felipe RS, Edwin MM Ortega, and Gauss M. Cordeiro."The generalized inverse Weibull distribution." Statistical Papers 52.3 (2011): 591-619.

Jenkinson, A. F. (1955). The frequency distribution of the annual maximum (or minimum) values of meteorological elements. Quarterly Journal of the Royal Meteorological Society, 81(348), 158-171.

Henningsen, A., &Toomet, O. (2011). maxLik: A package for maximum likelihood estimation in R. Computational Statistics, 26(3), 443-458.

Kersey, J. H., Weisdorf, D., Nesbit, M. E., LeBien, T. W., Woods, W. G., McGlave, P. B., Kim, T., Vallera, D. A., Goldman, A. I., Bostrom, B. and Ramsay, N. K. C. (1987). Comparison of autologous and autologous bone marrow transplantation for treatment of high-risk refractory acutelymphoblastic leukemia. New England Journal of Medicine 317, 461-467.

Krishna, E., Jose, K. K., Alice, T., & Ristić, M. M. (2013).The Marshall-Olkin Fréchet distribution.Communications in Statistics-Theory and Methods, 42(22), 4091-4107.

Louzada, F., Ramos, P. L.A New Long-Term Survival Distribution. Biostat Biometrics Open Acc J. 1(4): 555574.

Louzada, F, Ramos, P.L and Perdoná.G. C. (2016).Different estimation procedures for the parameters of the extended exponential geometric distribution for medical data. Computational and mathematical methods in medicine, v.2016, 1-12.

Maller, R. A., & Zhou, S. (1995). Testing for the presence of immune or cured individuals in censored survival data. Biometrics, 1197-1205.

Migon, H. S., Gamerman, D., & Louzada, F. (2014). Statistical inference: an integrated approach. CRC press.

Perdoná, G. C., & Louzada-Neto, F. (2011).A general hazard model for lifetime data in the presence of cure rate. Journal of Applied Statistics, 38(7), 1395-1405.

Ramos, P. L., Moala, F. A., & Achcar, J. A. (2014). Objective priors for estimation of extended exponential geometric distribution. Journal of Modern Applied Statistical Methods, 13(1), 226-243.

Ramos, P. L., Achcar, J.A., Moala, F. A., Ramos, E. & Louzada F., (2017) Bayesian analysis of the generalized gamma distribution using non-informative priors, Statistics 51(4), 824-843.

Roman, M., Louzada, F., Cancho, V. G., & Leite, J. G. (2012). A new long-term survival distribution for cancer data. Journal of Data Science, 10(2), 241-258.